\documentclass[reqno,a4paper]{amsart}
\usepackage[foot]{amsaddr}
\usepackage{graphicx,nicefrac}
\usepackage{algorithm,algorithmicx,algpseudocode}
\usepackage{a4wide}
\usepackage{ulem}
\algrenewcommand\algorithmicindent{1em}

\newcommand{\norm}[1]{\left\|#1\right\|}
\newcommand{\F}{\mathsf{F}}
\renewcommand{\L}{\mathcal{L}}
\renewcommand{\d}{\mathsf{d}}
\newcommand{\J}{\mathsf{J}}
\newcommand{\D}{\mathsf{D}}
\newcommand{\A}{\mathsf{A}}

\newcommand{\Id}{\,\mathsf{Id}}

\newtheorem{theorem}{Theorem}[section]
\theoremstyle{definition}
\newtheorem{example}[theorem]{Example}
\newtheorem{remark}[theorem]{Remark}
\title[Adaptive Newton-Type Schemes Based on Projections]{Adaptive Newton-Type Schemes Based on Projections}

\author[M.~Amrein]{Mario Amrein}
\address{Applied University of Zurich, CH-8401 Switzerland}
\email{mario.amrein@zhaw.ch}

\begin{document}
\normalem
\begin{abstract}
In this work we present and discuss a possible globalization concept for Newton-type methods. We consider nonlinear problems $f(x)=0$ in $\mathbb{R}^{n}$ using the concepts from ordinary differential equations as a basis for the proposed numerical solution procedure. Thus, the starting point of our approach is within the framework of solving ordinary differential equations numerically. Accordingly, we are able to reformulate general Newton-type iteration schemes using an adaptive step size control procedure. In doing so, we derive and discuss a discrete adaptive solution scheme thereby trying to mimic the underlying continuous problem numerically without losing the famous quadratic convergence regime of the classical Newton method in a vicinity of a regular solution. The derivation of the proposed adaptive iteration scheme relies on a simple orthogonal projection argument taking into account that, sufficiently close to regular solutions, the vector field corresponding to the Newton scheme is approximately linear. We test and exemplify our adaptive root-finding scheme using a few low-dimensional examples. Based on the presented examples, we finally show some performance data.
\end{abstract}



\maketitle

\section{Introduction}

In this note, we are interested in the problem: Find~$x_{\infty}\in \mathbb{R}^n$ such that 
\[f(x_{\infty})=0,\]
where~$f:\Omega \rightarrow \mathbb{R}^{n}$ denotes a possibly nonlinear function defined on the open subset~$\Omega \subset \mathbb{R}^{n}$. Of course this problem is one of the well known and possibly most addressed issues in numerical mathematics and has been studied by several authors in the past.
Here we study the problem of computing the roots of $f$ numerically. For~$x\in \Omega$ let the map~$x \mapsto \A(x) \in \mathbb{R}^{n \times n}$ be continuous. Next we set~$\F(x):=\A(x) f(x)$ and concentrate on the initial value problem

\begin{equation}
\label{eq:initial-value-problem}
\begin{cases} 
\begin{aligned}
 \dot{x}(t)&=\F(x(t)), \qquad && t\geq 0,\\
 x(0)&= x_0, \qquad && x_0 \in \Omega.
\end{aligned}
\end{cases}
\end{equation}

Assuming that a solution~$x(t)$ of~\eqref{eq:initial-value-problem} exists for all $t\geq 0$ with~$x(t)\in \Omega$, i.e.~$x_{\infty}:=\lim_{t\to \infty} {x(t)}\in \Omega$ and provided that~$\A(x_{\infty})f(x_{\infty})=0$ implies~$f(x_{\infty})=0$, we can try to follow the solution~$x(t)$ numerically to end up with an approximate root for~$f$. 
In actual computations however---apart from trivial problems---we can solve~\eqref{eq:initial-value-problem} only numerically. The simplest routine is given by the explicit---forward---Euler method. More precisely: For an initial value~$x_0\in \Omega$ a simple discrete version of the initial value problem~\eqref{eq:initial-value-problem} is given by
\begin{equation}
\label{eq:discrete-version}
x_{n+1}=x_{n}+t_{n}\F(x_n), \quad t_{n}\in (0,1], \quad n\geq 0.
\end{equation}
Obviously, depending on the non-linearity of~$\F$ and the choice of the initial value~$x_0$, such an iterative scheme is more or less meaningful for~$n \to \infty$. Indeed, supposing the limit for~$n\to \infty $ of the sequence~$(x_n)_{n\geq 0}$ generated by~\eqref{eq:discrete-version} exists, we end up with~$\F(x_n)\approx 0 $ for $n$ being sufficiently large.
Of course, we want to choose $\F$ in such a way that the iteration scheme in~\eqref{eq:discrete-version} is able to transport an initial value arbitrarily close to a root of~$\F$. 

For the remainder of this work, we assume that for all~$x_{n} $ generated by the iteration procedure from~\eqref{eq:discrete-version}, there exists a neighborhood of~$x_n$ such that the matrix~$ \A(x)$ is invertible.
Let us briefly address some different choices for~$\F$. A possible iteration scheme is based on~$\A(x):=~-\Id$ leading to a fixed point iteration which is also termed \emph{Picard iteration}. It is well known that under certain---quite strong---assumptions on~$f$ this scheme converges exponentially fast; see, e.g.,~\cite{Tao}. Another interesting choice for~$\F$ is given by~$\A(x):=\J_{f}(x)^{-1}$ leading to
\begin{equation}
\label{eq:choice1}
\F(x):=-\J_{f}(x)^{-1}f(x).
\end{equation}
Using this choice for~$\F$ in the iteration scheme~\eqref{eq:discrete-version} implies another well established iteration procedure called \emph{Newton's method} with damping. Here for~$x\in \Omega$ we denote by~$\J_{f}(x)$ the Jacobian of~$f$ at~$x$. Evidently this method requires reasonably strong assumptions with respect to the differentiability of $f$ as well as invertibility of the Jacobian~$\J_{f}(x_n)$ for all possible iterates~$x_n$ occurring during the iteration procedure. On the other hand---and on a \emph{local} level---, Newton's method with step size~$t_{n}\equiv 1$ is often celebrated for its super-exponential convergence regime `sufficiently' close to a regular root of~$f$.
Also well known are so called \emph{Newton-like methods} where the Jacobian~$\J_{f}(x)$ is replaced by a continuous approximation. A possible realization of such a method is given by setting~$\A(x):=-\J_{f}(x_0)^{-1}$, i.e., the initial derivative of~$f$ will be fixed through the whole iteration procedure. The iteration scheme~\eqref{eq:discrete-version} based on various choices for~$\F(x)$, where~$ \A(x)$ typically represents a (continuous) approximate of~$\J_{f}(x)^{-1}$ has been studied extensively by many authors in the recent past; see, e.g.,~\cite{5,epureanu:102,Ortega, Potschka}. Moreover, it is noteworthy that solving~\eqref{eq:initial-value-problem} with~$\F(x):=\J_{f}(x)^{-1}f(x)$ on the right is also known as the \emph{continuous Newton method}.
A pure analysis studying the long-term behavior of solutions for~\eqref{eq:initial-value-problem} which possibly lead to a solution of~$\F$ has also been studied in~\cite{1,2,3,smale2,Tanabe1}. Let us remark further that there is a wide research area where various methods are applied which are based on continuous Newton-type methods from~\eqref{eq:initial-value-problem} and its discrete analogue~\eqref{eq:discrete-version}.
The goal of the present work is not to give a complete summary of the wide-ranging theory and existing approaches for solving~\eqref{eq:initial-value-problem} within the context of a root finding procedure, but rather to illustrate some specific properties of vector fields $\F$ in order to understand the efficiency of the classical \emph{continuous Newton method} and thereby derive a simple---and efficient---adaptive numerical solution procedure for the numerical solution of the equation $f(x_{\infty})=0$. Although we only discuss the finite dimensional case, it is noteworthy that most of the following analysis extends without difficulty to the infinite dimensional case.

\subsubsection*{Notation.} In the main part of this paper, we suppose that---at least---there exists a zero~$x_{\infty}\in \Omega $ solving~$f(x_{\infty})=0$; here,~$\Omega $ denotes some open subset of the euclidean space~$\mathbb{R}^{n}$. In addition, for any two elements~$x,y \in  \mathbb{R}^{n} $ we signify by~$(x,y)$ the standard inner product of~$\mathbb{R}^{n}$ with the corresponding euclidean norm~$(x,x):=\norm{x}^{2}$. Moreover, for a given matrix~$\A \in \mathbb{R}^{n \times n}$ we denote by~$\norm{\A}$ the sup-norm induced by~$\norm{\cdot}$. We further denote by~$B_{R}(x)\subset \mathbb{R}^{n}$ the open ball with center at~$x$ and radius~$R>0$.
Finally, whenever the vector field~$f$ is differentiable, the derivative at a point~$x\in \Omega$ is written as~$\J_{f}(x)$, thereby referring to the Jacobian of~$f$ at~$x$.

\subsubsection*{Outline.}
This paper is organized as follows. In Section~\ref{sec:2} we first discuss the connection between the local and the global aspects of general Newton-type methods. More precisely, we interpret~$f$ as a vector field and focus on the local point of view, i.e, the case when an initial value $x_0 $ of the system~\eqref{eq:initial-value-problem} is `close' to a zero~$x_\infty$ of the vector field~$f$. Secondly, we consider the situation where initial guesses are no longer assumed to be `sufficiently close' to a zero~$x_\infty$ of~$f$. Based on the discussion within the local point of view, we transform the function~$f$ such that---at least on a local level---it is reasonable to expect convergence of our iteration scheme. In addition, we revisit the discretization of the initial value problem~\eqref{eq:initial-value-problem} in Section~\ref{sec:3} and define---based on the preceding results---an adaptive iteration scheme for the numerical solution of~\eqref{eq:initial-value-problem}. In Section~\ref{sec:4}, we present an algorithmic realization of the previous presented adaptive strategy. Finally, we give a series of low dimensional numerical experiments illustrating the performance of the adaptive strategy proposed in this work. Eventually, we summarize our findings in Section~\ref{sec:concl}. 
 

\section{Vector fields v.s. roots of a function $f$.}
\label{sec:2}
\subsection*{Local perspective:}
\begin{figure}[h!]
\includegraphics[width=0.4\textwidth]{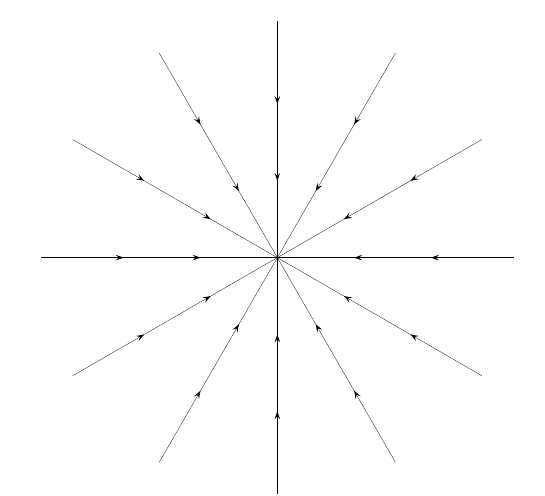}
\caption{The direction field associated with~$x\mapsto x_{\infty}-x$. Here, the center of the star signifies~$x_{\infty}$.}
\label{Fig:purely-linear}
\end{figure}

In this section, we start from a completely naive point of view by asking the following question: 
Is there a simple choice for the right hand side of~\eqref{eq:initial-value-problem} that can be used to transport an initial value~$x_{0} \in \Omega $ arbitrarily close to a root~$x_{\infty}\in \Omega $ of~$f$?
The answer to such a question typically depends on how close the initial guess~$x_0$ is chosen with respect to the zero~$x_{\infty}$. Indeed, if we assume that~$x_{0}$ is `sufficiently close' to~$x_{\infty}$ it would be preferable that such an initial guess~$x_0$ can be transported straightforwardly and arbitrarily close to the zero~$x_{\infty}$. However, let us remark on the following:

First of all and for the purpose of simplicity, we suppose that for~$x_{0}$ `sufficiently' close to~$x_{\infty}$ the function~$f$ is linear. More precisely, assume that the function~$f$ is given by~$f(x):=x_{\infty}-x$ (see also Figure~\ref{Fig:purely-linear}). Thus, if we set~$\A(x):=\Id $ on the right of~\eqref{eq:initial-value-problem}, the solution is given by~$x(t):=x_{_\infty}+(x_0-x_\infty)\mathrm{e}^{-t}$. Obviously, any initial guess~$x_0$ will be transported  arbitrarily close to the zero~$x_\infty$.
What can we learn from this favorable behavior of $x(t)$? On the one hand it would be preferable that an arbitrary vector field behaves like~$\F(x)=x_\infty-x$. In the nonlinear case and from a global perspective, i.e., whenever the initial guess~$x_0$ is far away from a zero~$x_{\infty}$, we would still like to establish a procedure which is able to transport the initial guess into a neighborhood of~$x_{\infty}$, where it is reasonable to assume the previous favorable behavior of the curve~$x(t)=x_{\infty}+(x_0-x_\infty)\mathrm{e}^{-t}$. So far, our discussion implies that on a local level we can typically expect to find a zero~$x_{\infty}$ whenever $\F(x)$ is close to~$x_{\infty}-x$. Let us therefore transform~$f$ in such a way that, at least on a local level,~$\F(x)=f(x)\approx x_\infty-x$ (see again Figure~\ref{Fig:purely-linear}) holds.

\subsection{Global perspective:}
As previously discussed, starting in~\eqref{eq:initial-value-problem} with an initial value~$x_{0}\in \Omega$, it would be preferable that the root~$x_\infty$ is attractive. More precisely, for the initial value~$x_0$ the corresponding solution~$x(t)$ should end at~$x_\infty$, i.e.,~$\lim_{t\to \infty}{x(t)}=x_\infty$ holds. Consequently, we would like to transform the vector field~$f$ in such a way that the new vector field---denoted by~$\F$---only has zeros which are at least `locally' attractive.
In other words, we want to transform~$f$ by~$\F(x):=\A(x)f(x)$ such that
\begin{equation}
\label{eq:w}
\F(x)\approx x_\infty-x, 
\end{equation}
holds true, especially whenever~$x$ is `close' to~$x_{\infty}$. A possible choice for~$\F$ that mimics the map~$x\mapsto x_{\infty}-x$ whenever~$x$ is close to the root~$x_\infty$ is given by 
\begin{equation}
\label{eq:Newton-flow}
\A(x):=-\J_{f}(x)^{-1}.
\end{equation}

Obviously, the price we have to pay for this choice is that~$f$ has to be differentiable with invertible Jacobian. Indeed, if~$f$ is twice differentiable with bounded second derivative, we observe that
\begin{equation}
\label{eq:Newton}
\begin{aligned}
\F(x)&=\F(x_\infty)+\D\F(x_{\infty})(x-x_\infty)+\mathcal{R}(x_{\infty},x-x_{\infty})\\
&=x_{\infty}-x+\mathcal{R}(x_{\infty},x-x_{\infty}),
\end{aligned}
\end{equation}
with~$ \norm{\mathcal{R}(x_{\infty},x-x_{\infty})}=\mathcal{O}(\norm{x-x_{\infty}}^2)$.

\begin{figure}
\includegraphics[width=0.5\textwidth]{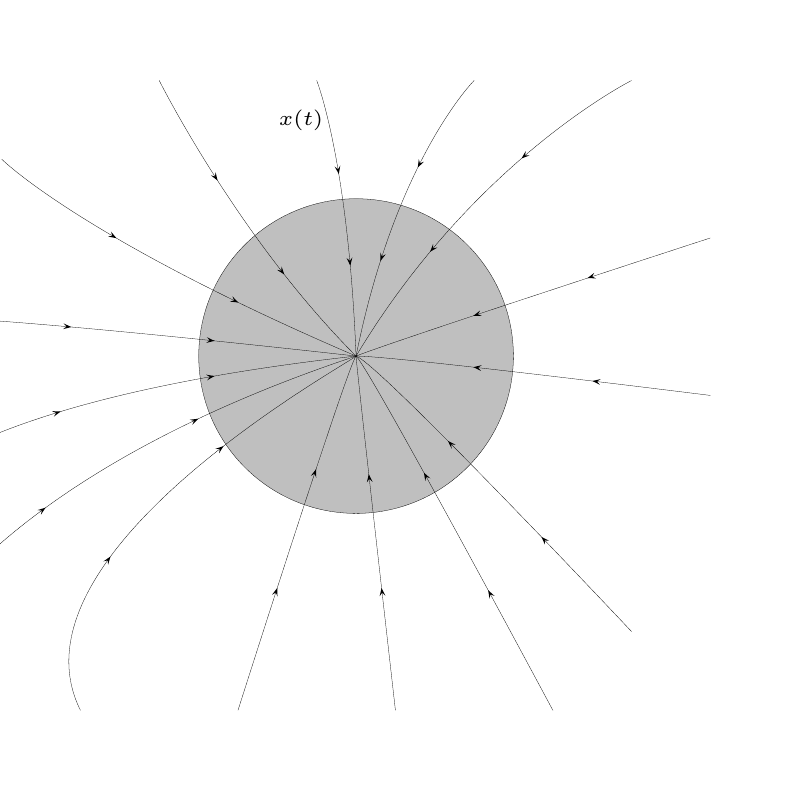}
\caption{An neighborhood of~$x_\infty$ where the map~$\F$ behaves like the affine linear map~$x\mapsto x_\infty-x$. Note that close to~$x_\infty$ the solutions are close to integral curves of the form~$x(t):=x_\infty+(x_0-x_\infty)\mathrm{e}^{-t}$ solving~\eqref{eq:initial-value-problem}.}
\label{Fig:Linear}
\end{figure}

Incidentally, it is well known that as long as the real parts of the eigenvalues of 
\[\D \F(x_{\infty})=\D[\A(x)f(x)]|_{x=x_{\infty}}=\A(x_\infty)\J_{f}(x_\infty),\] 
are negative, the zero~$x_\infty$ is locally attractive; see e.g.,~\cite{Königsberger}. As a result, if~$\A(x_\infty)$ is `sufficiently' close to the inverse of the Jacobian~$-\J_{f}(x_{\infty})$, the zero~$x_\infty$ might still be locally attractive. For example we can choose~$\F(x):=-\J_{f}(x_0)^{-1}f(x)$ in~\eqref{eq:discrete-version}. Generally speaking, whenever~$\A(x)$ is `sufficiently' close to the inverse of~$-\J_f(x)$ we still can hope that---especially on a local level---the iteration procedure~\eqref{eq:discrete-version} is well defined and possibly convergent, i.e.,~$x_{n}\to x_{\infty}$ for~$n \to \infty$.

Notice that whenever we can fix~$\A(x):=-\J_{f}(x)^{-1}$ the initial value problem given in~\eqref{eq:initial-value-problem} reads as follows:

\begin{equation}
\label{eq:cont-Newton}
\begin{cases} 
\begin{aligned}
 \dot{x}(t)&=-\J_{f}(x(t))^{-1}f(x(t)), \qquad && t\geq 0,\\
 x(0)&= x_0, \qquad && x_0 \in \mathbb{R}^{n}.
\end{aligned}
\end{cases}
\end{equation}

This initial value problem is also termed \emph{continuous Newton's method} and has been studied by several authors in the past; see, e.g.,~\cite{AmreinWihler:14, D04,epureanu:102,JLT07,neuberger,3,4,smale2,  Tanabe1, Tanabe2}.

Let us briefly show an important feature of the continuous Newton's method. Suppose that~$x(t)$ solves~\eqref{eq:cont-Newton}. Then it holds that
\[
\frac{\d }{\d t}f(x(t))=-f(x(t)),
\]
from which we deduce 
\[
f(x(t))=f(x_0)\mathrm{e}^{-t}.
\]

\section{Adaptivity based on a simple projection argument}
\label{sec:3}
In this section, we define an iteration scheme for the numerical solution of~\eqref{eq:initial-value-problem}. Based on the previous observations we further derive a computationally feasible adaptive step size control procedure. To this end, we assume that~$\F(x)= \A(x)f(x)$ is sufficiently smooth and that~$x_{\infty}$ is a regular root of~$f$, i.e.,~$\J_{f}(x_{\infty})^{-1}$ exists. Our analysis starts with a second order Taylor expansion of~$\F(x)$ around~$x_{\infty}$ given by
\begin{equation}
\label{eq:taylorexp}
\F(x)=\D \F(x_{\infty})(x-x_{\infty})+\mathcal{R}_{x_{\infty}}(x-x_{\infty}), \quad \norm{\mathcal{R}_{x_{\infty}}(x-x_{\infty})}=\mathcal{O}(\norm{x-x_{\infty}}^2).
\end{equation}

Next we recall that whenever we are able to choose 
$\A(x):=-\J_{f}(x)^{-1}$ there holds 
\[
\F(x)=\L(x)+\mathcal{R}_{x_{\infty}}(x-x_{\infty}),
\]
where $ \L(x):=x_{\infty}-x$.

For~$x,y \in \mathbb{R}^{n}$ we consider the orthogonal projection of~$x$ onto~$y$ given by

\begin{equation}
\label{eq:proj}
\text{proj}_{y}(x):=\frac{(x,y)}{\norm{y}^2}\cdot y.
\end{equation}

Furthermore, we use the fact that in a neighborhood of a regular root~$x_{\infty}$ it holds that~$\F(x)\approx \L(x)$. Based on these observations, we use the orthogonal projection of~$\F(x)$ onto~$\L(x)$. In particular for~$x\neq x_{\infty}$~\eqref{eq:taylorexp} delivers
\[
\text{proj}_{\L(x)}(\F(x))=\frac{(\F(x),\L(x))}{\norm{\L(x)}^2}\cdot \L(x)=\L(x)+\frac{(\mathcal{R}_{x_{\infty}}(x-x_\infty),\L(x))}{\norm{\L(x)}^2}\cdot \L(x).
\]
Evidentially in case of $ \F(x)=\L(x)$ there holds $\text{proj}_{\L(x)}=\Id$.

For the remainder of this section, we assume that for an initial guess~$x_{0}\in \Omega $ there exists a solution~$x(t)$ for the initial value problem from~\eqref{eq:initial-value-problem} such that~$\lim_{t\to \infty}{x(t)}=x_{\infty}$ solves~$f(x_{\infty})=0$. Moreover, we assume that there exists an open neighborhood~$B_{R}(x_{0})\subset \Omega$ of $x_{0}$ such that for all~$x\in B_{R}(x_{0})$ there exists a solution~$x(t)$ of~\eqref{eq:initial-value-problem} starting in~$x\in B_{R}(x_{0})$ with~$\lim_{t\to \infty}{x(t)}=x_{\infty}$.
Thus, for~$t>0 $ sufficiently small, we can assume that 
\begin{equation}
\label{eq:x1x2}
x_{1}:=x_0+t\F(x_0) \quad \text{and} \quad x_{2}:=x_1+t\F(x_1).
\end{equation}
are elements of~$B_{R}(x_{0})$.

\begin{figure}
\includegraphics[width=0.6\textwidth]{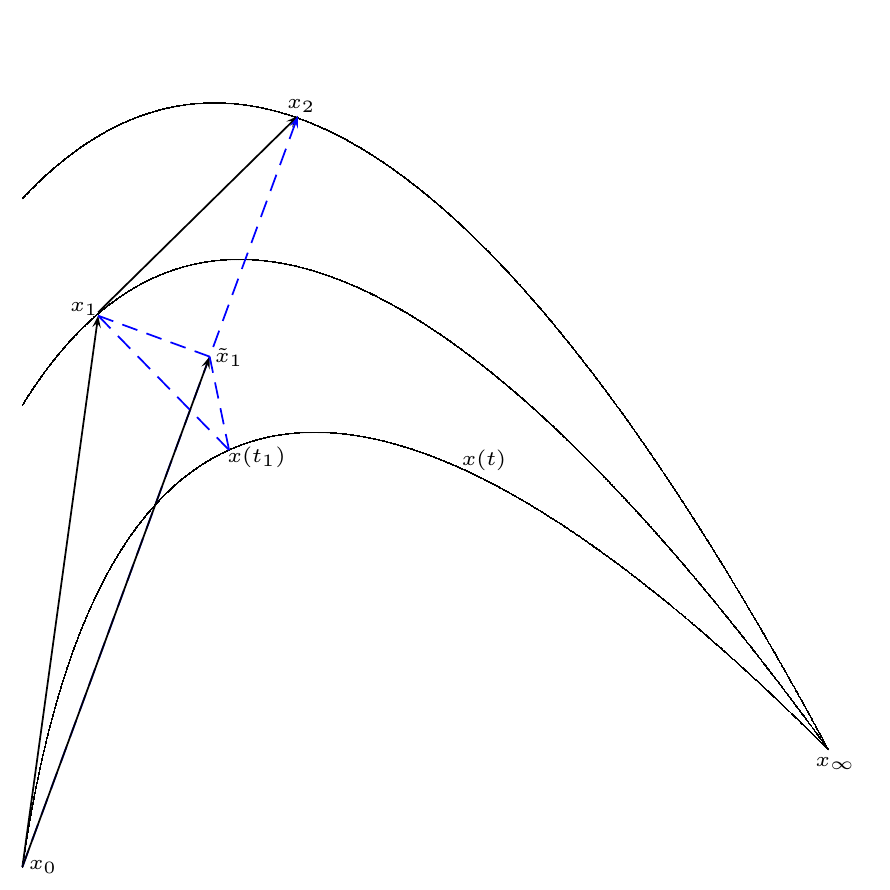}
\caption{The projection of $x_{1}-x_{0}=t_{1}\F(x_0)$ onto $t_1v=x_{2}-x_{0}$ after a time step $t=t_{1}$.}
\label{projection}
\end{figure}

We now use~$ \F(x_{0}) $ and $\F(x_{1})$ and set~$v:=\F(x_0)+\F(x_1)$. Note that for~$t=1$ and $x_{2}$ `close' to~$x_{\infty}$ there holds~$v=x_2-x_0\approx x_{\infty}-x_{0}=\L(x_0)$. 
Next we define our effectively computed iterate
\begin{equation}
\label{eq:eff}
\tilde{x}_1:=x_0+t\text{proj}_{v}(\F(x_0)).
\end{equation}
The situation is depicted in Figure~\ref{projection}. Note that~$\F(x_1)\approx \F(x_0)$ implies
\[
\text{proj}_{v}(\F(x_0))=\frac{(v,\F(x_0))}{\norm{v}^2}v\approx \F(x_0).
\]

Let us focus on the distance between the exact solution~$x(t)$ and its approximate~$\tilde{x}_{1}$ at $t>0$. In doing so we revisit the proposed approach from~\cite[\S2.3]{AmreinWihler:15}.

\subsection*{Error Analysis}
First we consider the Taylor expansion of~$x(t)$ around~$t_0=0$:
\begin{equation}
\label{eq:1}
\begin{aligned}
x(t)&=x_{0}+t\dot{x}(0)+t^2\frac{\ddot{x}(0)}{2}+\mathcal{R}_{x}(t)\\
&=x_0+t\F(x_0)+t^2\frac{\ddot{x}(0)}{2}+\mathcal{R}_{x}(t)\\
&=x_1+t^2\frac{\ddot{x}(0)}{2}+\mathcal{R}_{x}(t), \quad \text{with} \quad \norm{\mathcal{R}_{x}(t)} = \mathcal{O}(t^3).
\end{aligned}
\end{equation}
Moreover, we will take a look at the expansion of~$\F(x)$ around~$x_1$ given by
\[
\F(x_1)=\F(x_0)+t\D\F(x_0)\F(x_0)+\mathcal{R}_{\F}(t\F(x_0)) \quad \text{with} \quad \norm{\mathcal{R}_{\F}(t\F(x_0))}=\mathcal{O}(t^2\norm{\F(x_0)}^2).
\]

We see that there holds
\[
\lim_{t \searrow 0}{\frac{\F(x_1)-\F(x_0)}{t}}=\D\F(x_0)\F(x_0)=\frac{\d}{\d t}\F(x(t))|_{t=0}=\ddot{x}(0),
\]
and therefore 
\begin{equation}
\label{eq:2}
\F(x_1)-\F(x_0)=t\ddot{x}(0)+\mathcal{R}_{\F}(t\F(x_0)).
\end{equation}

Next we employ~\eqref{eq:1} and~\eqref{eq:2} in order to end up with
\begin{equation}
\label{eq:3}
x(t)-x_1=t^2\frac{\ddot{x}(0)}{2}+\mathcal{R}_x(t)=\frac{t}{2}(\F(x_1)-\F(x_0))+t\mathcal{R}_{\F}(t\F(x_0))+\mathcal{R}_{x}(t).
\end{equation}

\begin{remark}
Note that~\eqref{eq:3} can serve as an error indicator for the iteration~\eqref{eq:discrete-version} (see~\cite[\S2.3]{AmreinWihler:15} or~\cite[\S2.2]{AmreinWihlerMelenk} for further details).
\end{remark}

Now we consider the difference~$x(t)-\tilde{x}_{1}$ using~\eqref{eq:3}:

\begin{align*}
x(t)-\tilde{x}_1&=x(t)-x_1+x_1-\tilde{x}_1\\
&=\frac{t}{2}(\F(x_1)-\F(x_0))+t(\F(x_0)-\text{proj}_{v}(\F(x_0)))+t\mathcal{R}_{\F}(t\F(x_0))+\mathcal{R}_{x}(t)\\
&=t\left(\frac{v}{2}-\text{proj}_{v}(\F(x_0))\right)+t\mathcal{R}_{\F}(t\F(x_0))+\mathcal{R}_{x}(t).
\end{align*}

Hence if we set~$\gamma(x_0,x_1):=\norm{\nicefrac{v}{2}-\text{proj}_{v}(\F(x_0))}$ we get
\begin{equation}
\label{eq:basic}
\norm{x(t)-\tilde{x}_1}\leq t\gamma(x_0,x_1)+\mathcal{O}(t^3).
\end{equation}

We see that by neglecting the term~$\mathcal{O}(t^3)$, the expression~$t\gamma(x_0,x_1)$ can be used as an error indicator in each iteration step. Furthermore, for~$\F(x_1)=\F(x_0)$ there holds~$\gamma(x_0,x_1)=0$.

Thence, fixing a tolerance~$\tau >0 $ such that 
\begin{equation}
\label{eq:equation}
\tau = t\gamma(x_0,x_1),
\end{equation}
 
motivates an adaptive step size control procedure for the proposed iteration scheme given in~\eqref{eq:eff} that will be discussed and tested in the next section.

\subsection{Adaptive strategy}
We now propose a procedure that realizes an adaptive strategy based on the previous observations. The individual computational steps are summarized in Algorithm~\ref{al:full}.

\begin{algorithm}
\caption{Adaptive Newton-like method:}
\label{al:full}
\begin{algorithmic}[1]
\State{\textbf{Input}:}
\begin{enumerate} 
\item[$\bullet$]
initial value $x_0 \in \Omega$,
\item[$\bullet$]
lower bound for the step size $t_{\text{lower}}>0$, 
\item[$\bullet$]
error tolerance $\tau>0$ and $\varepsilon>0$ respectively.
\end{enumerate}
\State{$\F(x_0)\leftarrow \A(x_0)f(x_0)$}
\State{$t\leftarrow \min\left(1,\sqrt{\frac{2\tau}{\norm{\F(x_0)}}}\right)$}
\For{$k=1,2,\dots$}
\If{$\norm{\F(x_0)}\leq \varepsilon$}
\State{\Return $x_{\infty} \leftarrow x_0$}
\Else
\Loop \Comment{start the adaptive step size control}
\If{$t<t_{\text{lower}}$}
\State{\textbf{stop the iteration procedure}}
\EndIf
\State{$x_1\leftarrow x_0+t\F(x_0)$}
\State{$\F(x_1)\leftarrow \A(x_1)f(x_1)$}
\State{$v\leftarrow \F(x_1)+\F(x_0)$}
\State{$\text{proj}_{v}(\F(x_0))\leftarrow \frac{(v,\F(x_0))}{\norm{v}^2}v$}
\State{$\gamma(x_0,x_1)\leftarrow \norm{\nicefrac{v}{2}-\text{proj}_{v}(\F(x_0))}$}
\If{$t\gamma(x_0,x_1)\leq\tau $ } 
\State{\textbf{break the loop}}
\Else
\State{{$t\leftarrow \mathsf{R}(t)$}{\Comment reduce the step size}}
\EndIf
\EndLoop
\State{{$x_0\leftarrow x_0+t\text{proj}_{v}(\F(x_0))$}{\Comment perform a step}}
\State{{$\F(x_0)\leftarrow \A(x_0)f(x_0)$}{\Comment update the direction}}
\State{{$t\leftarrow \min{\left(1,\frac{\tau}{\gamma(x_0,x_1)}\right)}$}{\Comment predict the step size}}
\EndIf
\EndFor
\end{algorithmic}
\end{algorithm}

\begin{remark}
By~$\mathsf{R}(t)$ we signify a procedure that reduces the current step size such that~$0<\mathsf{R}(t)<t$.
Let us also briefly address a possible and reasonable choice for the initial step size~$t_{\text{init}}$ in Algorithm~\ref{al:full}. 
The following---detailed---reasoning can also be found in~\cite[\S2]{AmreinWihler:15}.

If we start our procedure with a regular initial value~$x_{0}\in\Omega$ such that~$\F(x)=-\J_{f}(x)^{-1}f(x)$ is Lipschitz continuous in a neighborhood of~$x_0$, then there exists a local---unique---solution for~\eqref{eq:initial-value-problem}, i.e., there exists~$T>0$ with 
\[
\dot{x}(t)=-\J_{f}(x(t))^{-1}f(x(t)),
\]
on~$t\in [0,T)$. Consequently, there holds~$f(x(t))=f(x_0)\mathrm{e}^{-t}$. A second order Taylor expansion reveals (see~\cite[\S2.2]{AmreinWihler:15} or~\cite[\S2.2]{AmreinWihler:14} for details) 
\begin{equation}
\label{eq:second}
x(t)\approx x_{0}+\dot{x}(0)t+t^2\xi=x_{0}+t\F(x_0)+t^2\xi,
\end{equation}
with~$\xi \in \mathbb{R}^{n}$ to be determined. Moreover we use a second order Taylor expansion for~$f$ and compute
\begin{equation}
\label{eq:init}
f(x_0)\mathrm{e}^{-t}=f(x(t))\approx f(x_{0}+\dot{x}(0)t+t^2\xi)
=f(x_0)-tf(x_0)+\J_{f}(x_0)t^2\xi.
\end{equation}

We finally use~$\mathrm{e}^{-t}\approx 1-t+\frac{t^2}{2}$ in 
\[
f(x_0)\mathrm{e}^{-t}\approx f(x_0)-tf(x_0)+\J_{f}(x_0)t^2\xi
\]
in order to end up with
\[
\xi\approx \frac{1}{2}\J_{f}(x_0)^{-1}f(x_0).
\]
Combining this with~\eqref{eq:second} yields
\[
x(t)\approx x_0+t\F(x_0)+\frac{1}{2}t^2\J_{f}(x_0)^{-1}f(x_0).
\]
Note that~$x_{1}=x_0+t\F(x_0)$. Thus after a first step~$t=t_{\text{init}}>0$ we get 
\[
\norm{x(t)-x_1}\approx \frac{1}{2}t^2\norm{\J_{f}(x_0)^{-1}f(x_0)}.
\]
Thus for the given error tolerance~$\tau >0 $, we set
\[
t_{\text{init}}=\sqrt{\frac{2\tau}{\norm{\J_{f}(x_0)^{-1}f(x_0)}}},
\]
i.e., we arrive at~$\norm{x(t_{\text{init}})-x_1}\approx \tau $.
\end{remark}

\begin{remark}
In Algorithm~\ref{al:full} the minimum in Step 3 and 25 respectively is chosen such that~$t=1$ whenever possible, in particular, whenever the iterates are close to the zeros~$x_\infty$. This will retain the famous quadratic convergence property of the classical Newton scheme (provided that the zero~$x_{\infty}$ is simple).
\end{remark}

\section{Numerical Experiments}
\label{sec:4}

The purpose of this section is to illustrate the performance of Algorithm~\ref{al:full} by means of two examples. In particular, we consider three algebraic systems. The first one is a polynomial equation on~$\mathbb{C} $ (identified with~$\mathbb{R}^{2}$) with three separated zeros, and the second example is a challenging benchmark problem in~$\mathbb{R}^{2}$. Finally we consider a problem in~$\mathbb{R}^{2}$ with exactly one zero in order to highlight the fact that---in certain situations---the classical Newton method is able to find a numerical solution whereas the proposed adaptive scheme is not convergent. 

For all presented examples, we set in Algorithm~\ref{al:full}

\[
t_{\text{init}}=\min\left(\sqrt{\frac{2\tau}{\norm{\J_{f}(x_0)^{-1}f(x_0)}}},1\right).
\]

In Algorithm \ref{al:full}, the lower bound for the step size in Step 9 is set to~$t_{\text{lower}}=10^{-9}$ and for the error tolerance in Step 5 we use~$\varepsilon=10^{-8}$. Moreover, for the possible reduction procedure~$t \leftarrow \mathsf{R}(t)$ in Step 20 we simply use~$\mathsf{R}(t):=\frac{t}{2}$. Let us further point out that there are more sophisticated strategies for the reduction process of the time step~ $t$ (see also~\cite[\S10]{Potschka}). Finally, we set the maximal number of iterations~$n_{\max}$ to 100.

\begin{example}
\label{ex:1}
We consider the function
\[
f:\mathbb{C}\rightarrow \mathbb{C}, \quad z \mapsto f(z):=z^3-1.
\]
Here, we identify~$f$ in its real form in~$\mathbb{R}^{2}$, i.e., we separate the real and imaginary parts. The three zeros are given by
\[
Z_{f}=\{(1,0),(-\nicefrac{1}{2},\nicefrac{\sqrt{3}}{2}),(-\nicefrac{1}{2},-\nicefrac{\sqrt{3}}{2})\} \subset \mathbb{C}.
\]

Note that~$\J_{f}$ is singular at~$(0,0)$. Thus if we apply the classical Newton method with~$\F(x)=-\J_{f}(x)^{-1}f(x)$ in~\eqref{eq:discrete-version} the iterates close to~$(0,0)$ causes large updates in the iteration procedure. More precisely, the application of~$\F(x)=-\J_{f}(x)^{-1}f(x)$ is a potential source for chaos near~$(0,0)$. Before we discuss our numerical experiments, let us first consider the vector fields generated by the continuous problem~\eqref{eq:initial-value-problem}. In Figure~\ref{fig:flows1}, we depict the direction fields corresponding to~$\F(x)=f(x)$ (left) and~$ \F(x)=-\J_{f}(x)^{-1}f(x)$ (right). We clearly see that~$(0,1)\in Z_{f}$ is repulsive for~$\F(x)=f(x)$. Moreover, the zeroes~${(-\nicefrac{1}{2},\nicefrac{\sqrt{3}}{2}),(-\nicefrac{1}{2},-\nicefrac{\sqrt{3}}{2})}\in Z_{f}$ of~$\F(x)=f(x)$ show a curl. If we now consider~$ \F(x)=-\J_{f}(x)^{-1}f(x)$ the situation is completely different: All zeros are obviously attractive. In this example, we further observe that the vector direction field is divided into three different sectors, each containing exactly one element of~$Z_{f}$. 

Next we visualize the domains of attraction of different Newton-type schemes. In doing so, we compute the zeros of~$f$ by sampling initial values on a~$500\times 500 $ grid in the domain~$[-3,3]^2$ (equally spaced). In Figure~\ref{fig:frac12}, we show the fractal generated by the traditional Newton method with step size~$t\equiv 1$ (left) as well as the corresponding plot for the adaptive Newton-type scheme with the proposed variable step size~$t$ (right). It is noteworthy that the chaotic behavior caused by the singularities of~$\J_{f}$ is clearly tamed by the adaptive procedure. Here, we set~$\tau = 0.01$ in Algorithm~\ref{al:full}. 

\begin{figure}
\includegraphics[width=0.45\textwidth]{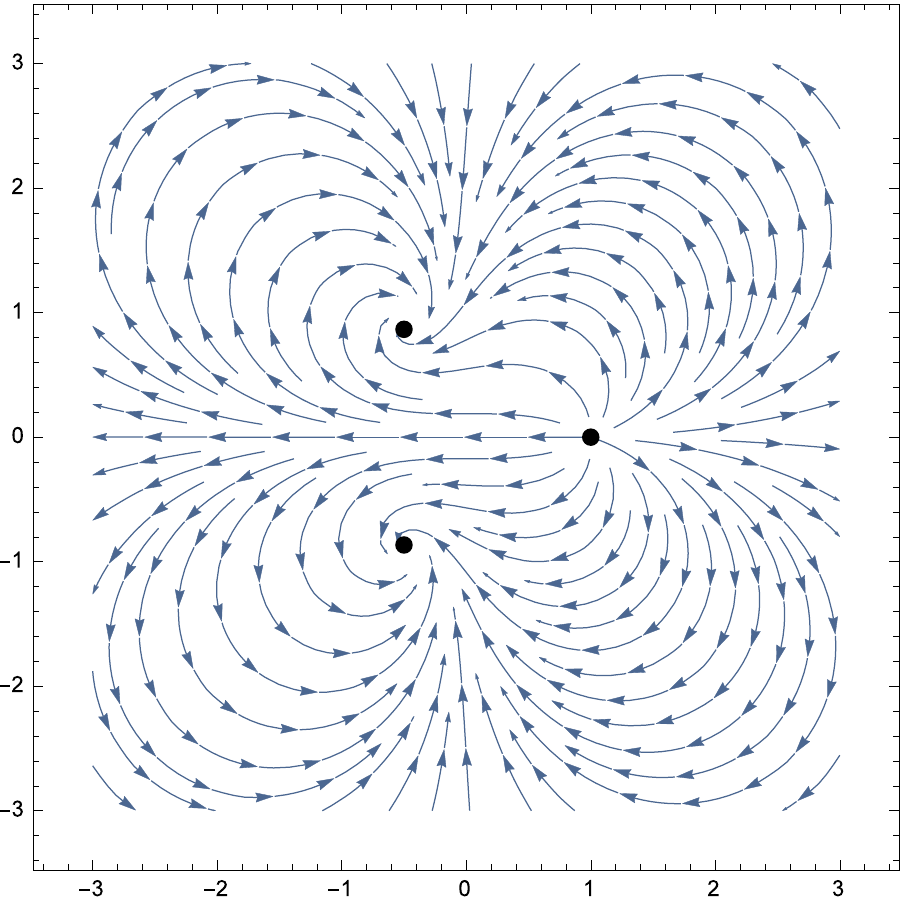}
\hfill
\includegraphics[width=0.45\textwidth]{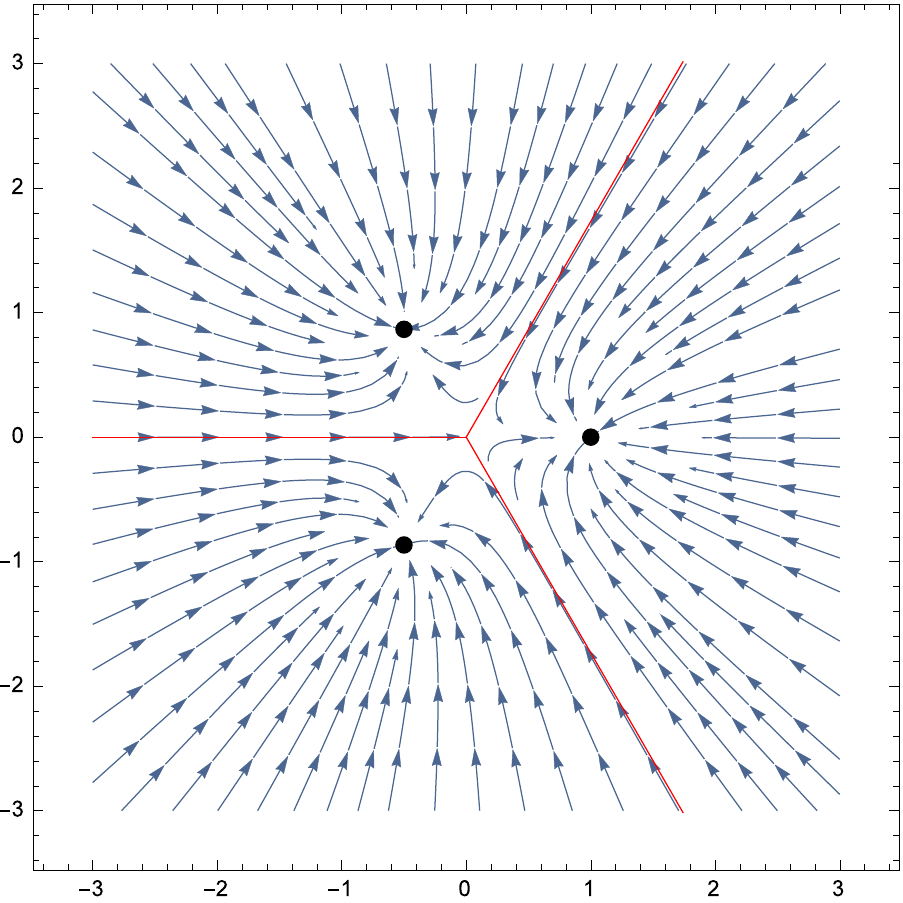}
\caption{Example~\ref{ex:1}: The direction fields corresponding to~$f(z)=z^3-1$ (left) and to the transformed~$\F(z)=-\J_{f}(z)^{-1}\cdot f(z)$ (right).}
\label{fig:flows1}
\end{figure}

\begin{figure}
\includegraphics[width=0.452\textwidth]{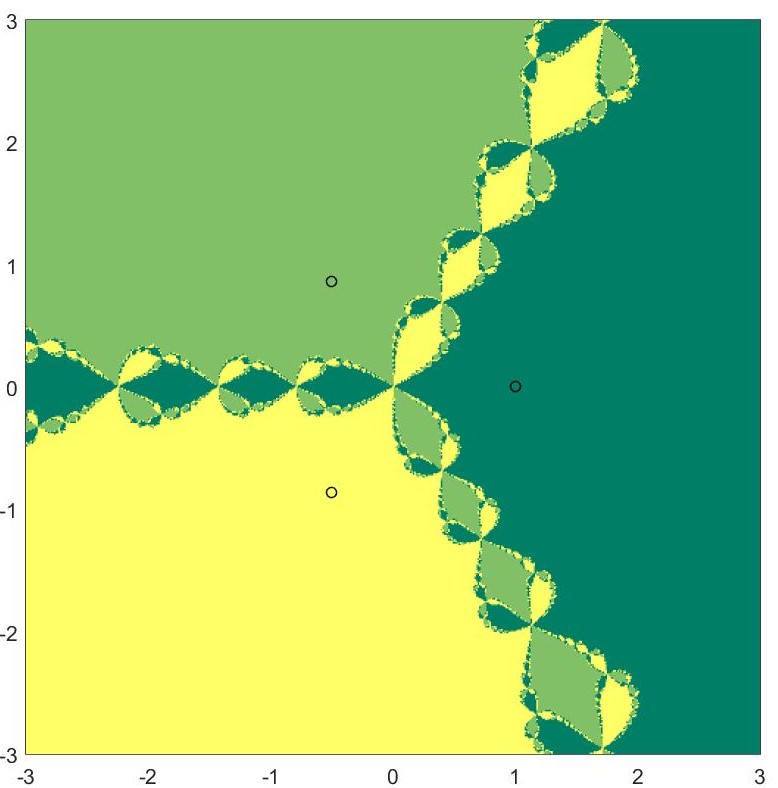}
\hfill
\includegraphics[width=0.447\textwidth]{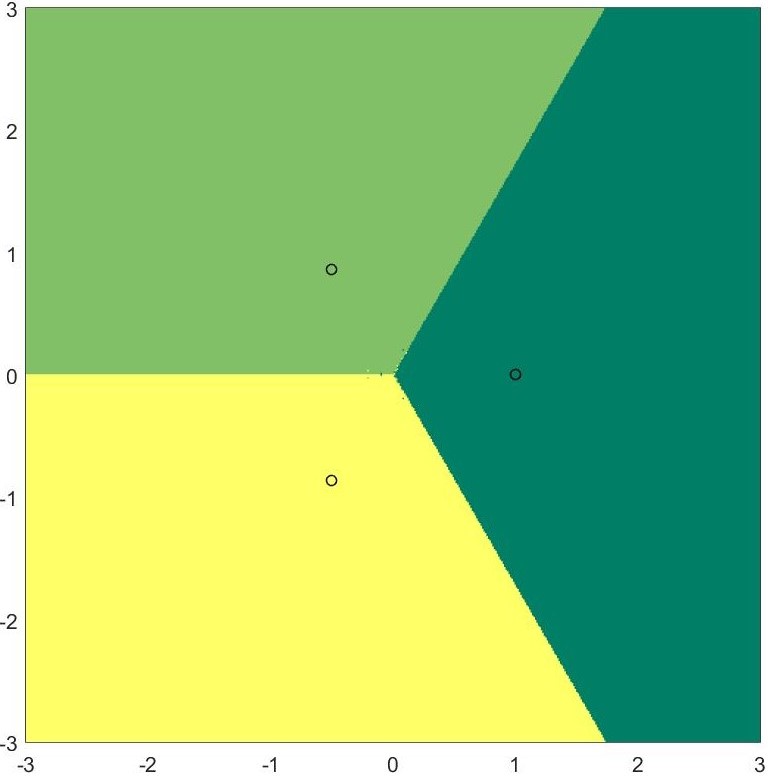}
\caption{The basins of attraction for Example~\ref{ex:1} by the Newton method. On the left without step size control (i.e.,~$t\equiv1$) and on the right with step size control~($\tau = 0.01$). Three different colors distinguish the three basins of attraction associated with the three solutions (each of them is marked by a small circle).}
\label{fig:frac12}
\end{figure}

\end{example}

\begin{example}
\label{ex:2}

The second test example is a~$2\times 2$ algebraic system from~\cite{5} defined as
\begin{equation}
\label{eq:benchmark}
f:[-1.5,1.5]^{2}\rightarrow \mathbb{R}^{2}, \quad f(x,y):=
\begin{pmatrix}
\text{exp}(x^2+y^2)-3 \\
x+y-\sin(3(x+y))
\end{pmatrix}.
\end{equation}

Firstly we notice that the singular set for~$\J_{f}$ is given by 
\[
\{y=x\}, \quad \text{and} \quad \left\{y=-x\pm \frac{1}{3}\arccos\left(\frac{1}{3}\right)\pm \frac{2}{3}\pi k, k \in \mathbb{N}_{\geq 0}\right\}.
\]

In Figure~\ref{fig:flows2}, we depict again the direction field associated to~$\F(x)=f(x)$ (left) and~$\F(x)=-\J_{f}(x)^{-1}f(x)$ (right).
If we use~$\F(x)=-\J_{f}(x)^{-1}f(x)$, we clearly see that six different zeros of~$f$ are all locally attractive. The right lines in Figure~\ref{fig:flows2} (right) indicate the singular set of~$\J_{f}$. In Figure~\ref{fig:34}, we show the domain of attraction. We clearly see that the proposed adaptive scheme in Algorithm~\ref{al:full} is able to tame the chaotic behavior of the classical Newton iteration. Let us further point out the following---important---fact: 

Suppose we are given an initial value~$x_{0}$ for the continuous problem~\eqref{eq:initial-value-problem} which is located in the subdomain of~$[-1.5,1.5]^{2}$ where no root of~$f$ is located (see the upper right and the bottom left part of the domain~$[-1.5,1.5]^{2}$ in Figure~\ref{fig:flows2} right). The trajectories corresponding to such initial guesses end at the singular set of~$\J_{f}$. The situation is different in the discrete case. Indeed, if we start the Newton-type iteration in~\eqref{eq:discrete-version} on the subdomain where no zero of~$f$ is located, the discrete iteration is potentially able to cross the singular set. In addition, if we set~$\tau \ll 1$ the discrete iteration~\eqref{eq:discrete-version} is close to its continuous analogue~\eqref{eq:initial-value-problem}. Therefore a certain amount of chaos may enlarge the domain of convergence. This is particularly important when no a priori information on the location of the zeros of~$f$ is available. We depict this situation in the Figures~\ref{fig:34}. Here we sample~$250\times 250$ equally spaced initial guesses on the domain~$[-1.5,1.5]^2$. The dark blue shaded part indicates the domain where the iteration fails to converge. Note that the proposed step size control is able to reduce the chaotic behavior of the classical Newton method. Moreover, the domain of convergence is again considerably enlarged by the adaptive iteration scheme.

\begin{figure}
\includegraphics[width=0.45\textwidth]{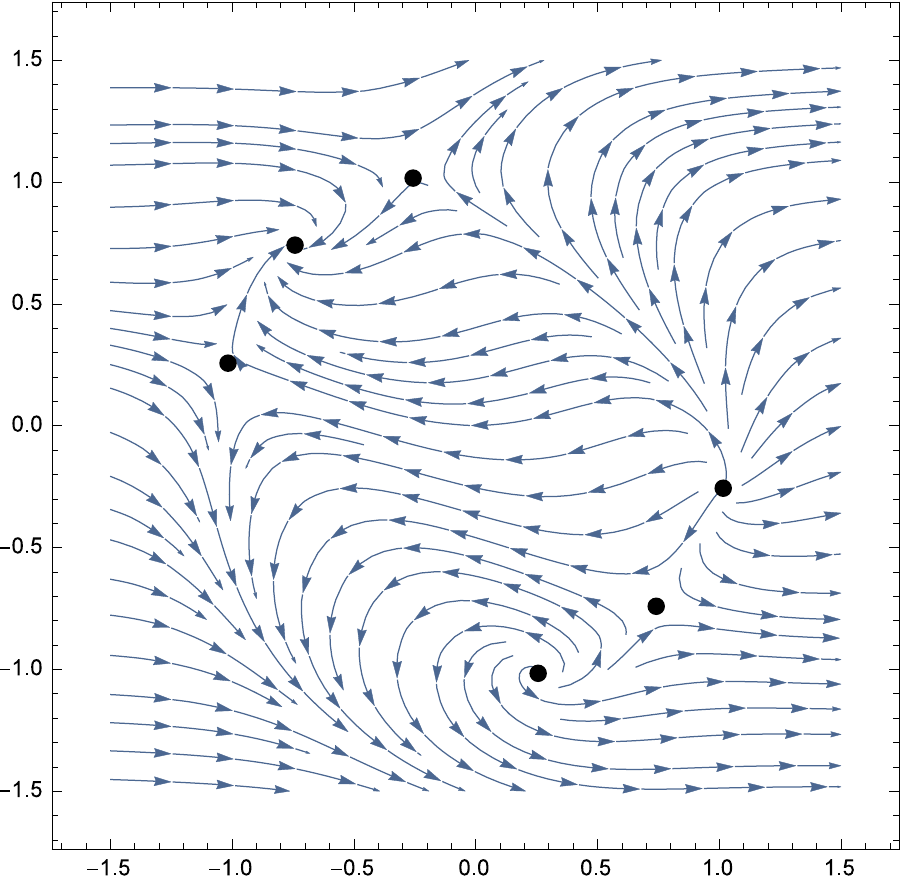}
\hfill
\includegraphics[width=0.45\textwidth]{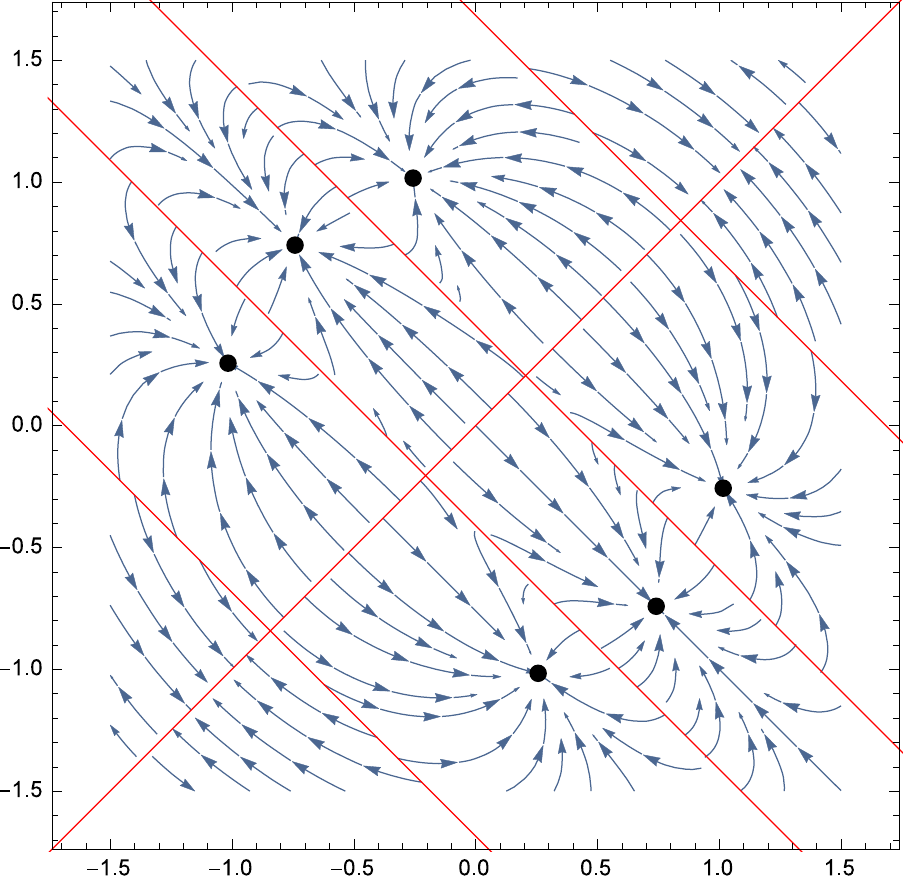}
\caption{Example~\ref{ex:2}: The direction fields corresponding to Example~\ref{ex:2}. On the left for~$\F(x)=f(x)$ and to the right for the  transformed vector field~$\F(x)=-\J_{f}(x)^{-1} f(x)$.}
\label{fig:flows2}
\end{figure}

\begin{figure}
\includegraphics[width=0.45\textwidth]{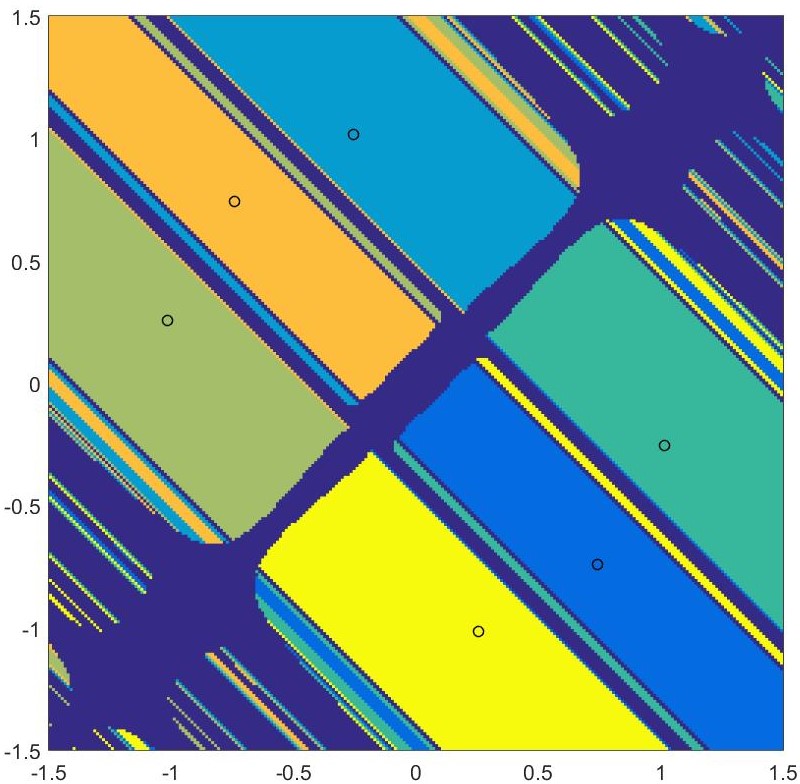}
\hfill
\includegraphics[width=0.448\textwidth]{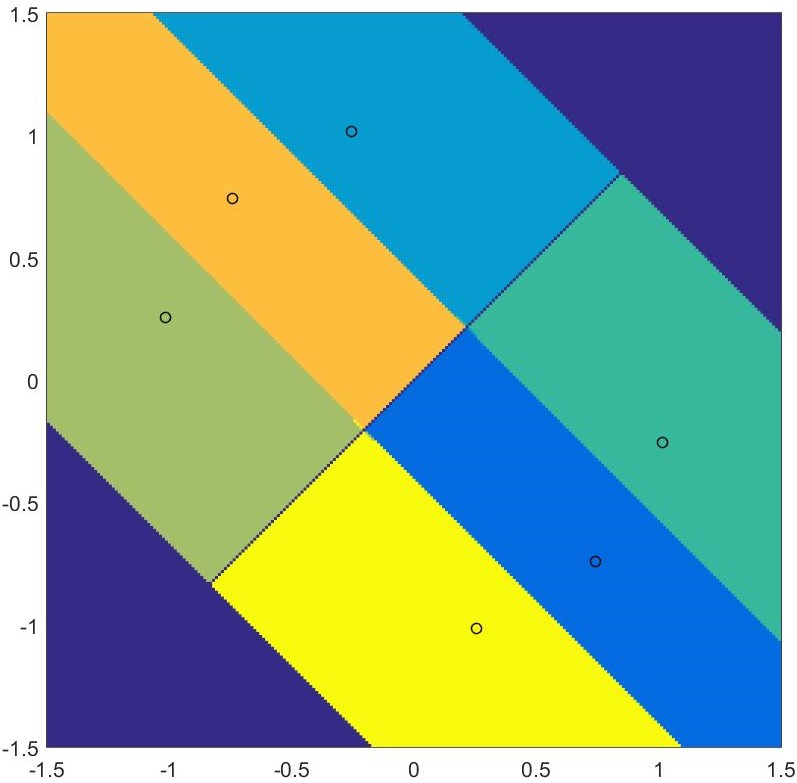}
\caption{The basins of attraction for Example~\ref{ex:2} by the Newton method. On the left without step size control (i.e,~$t=1$) and on the right with step size control~($\tau = 0.01$). Six different colors distinguish the six basins of attraction associated with the six solutions (each of them is marked by a small circle). Note that the dark-blue shaded domain indicates the domain, where the iteration procedure~\ref{al:full} fails to convergence (within the maximal number of iterations which is set here to~$n_{\text{max}}=100)$.}
\label{fig:34}
\end{figure}

\end{example}

\begin{example}
\label{ex:3}

We finally consider the algebraic $2\times 2$ system from~\cite{ScWi11} given by
\begin{equation}
\label{schnee}
f:[-10,10]^2\rightarrow \mathbb{R}^2, 
\quad f(x,y):=
\begin{pmatrix}
-x^2+y+3 \\
-xy-x+4
\end{pmatrix}.
\end{equation}

There exists a unique zero for~$f$ given by~$(2,1)$. This zero is an attractive fixed point for the system~\eqref{eq:initial-value-problem} with~$\F(x)=f(x)$ as well as~$\F(x)=-\J_{f}(x)^{-1}f(x)$. The associated direction fields are depicted in Figure~\ref{fig:fluss56}. Close to the zero~$(2,1)$ we observe a curl in case of~$\F(x)=f(x)$. However, if we instead use~$\F(x)=-\J_{f}(x)^{-1}f(x)$, the curl is removed and the direction field points directly to~$(2,1)$. In Figure~\ref{fig:frac56}, we show the attractors of~$(2,1)$ for the classical Newton method (left) and for the proposed adaptive strategy with~$\tau = 0.01$ (right). These pictures are based on sampling~$10^{6}$ starting values in the domain~$[-10,10]^2$. The right and yellow shaded part signifies the attractor for~$(2,1)$. Again we notice that the classical Newton method with step size~$t\equiv 1$ produces chaos. In the adaptive case the situation is different. We clearly see that adaptivity again is able to reduce the chaos and unstable behavior of the classical Newton method. Referring to the previous Example~\ref{ex:2}, it is noteworthy that in Example~\ref{ex:3} the domain of convergence in the adaptive case is comparable to the case of~$t\equiv 1$, i.e., the classical Newton method. 

\begin{figure}
\includegraphics[width=0.45\textwidth]{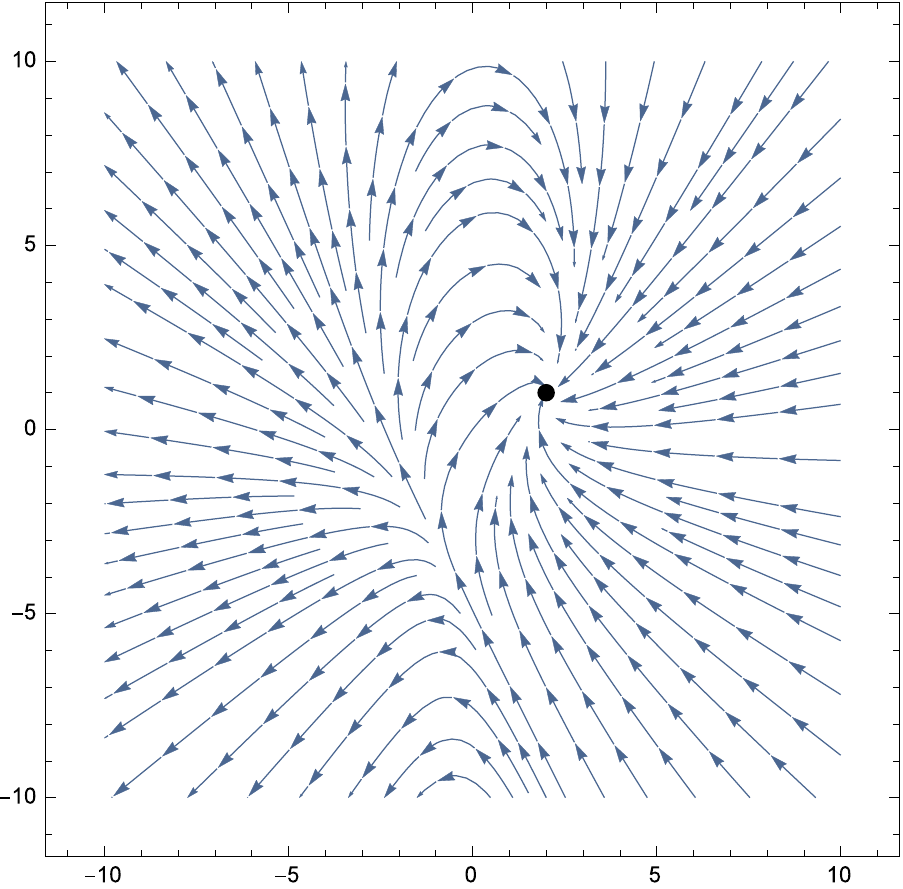}
\hfill
\includegraphics[width=0.45\textwidth]{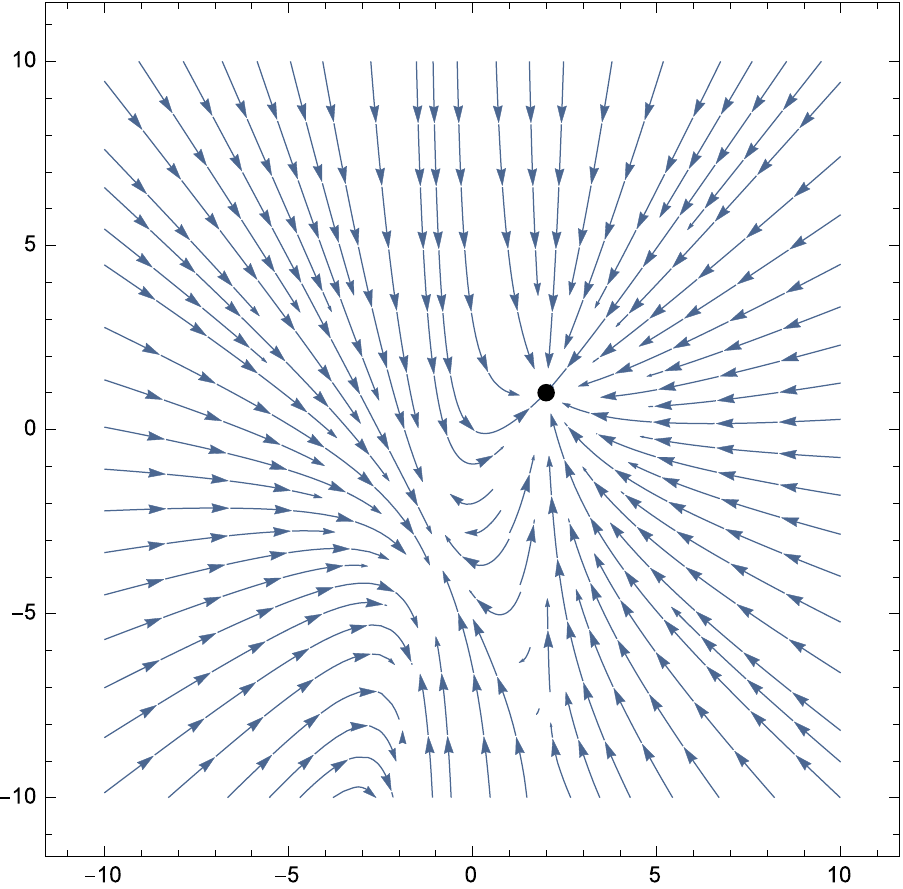}
\caption{The direction fields corresponding to Example~\ref{ex:3}. On the left for~$\F(x)=f(x)$ and on the right for the transformed vector field~$\F(x)=-\J_{f}(x)^{-1}f(x)$. We clearly see that the transformed field removes the curl which we obtain by simply applying~$\F(x)=f(x)$.}
\label{fig:fluss56}
\end{figure}

\begin{figure}
\includegraphics[width=0.45\textwidth]{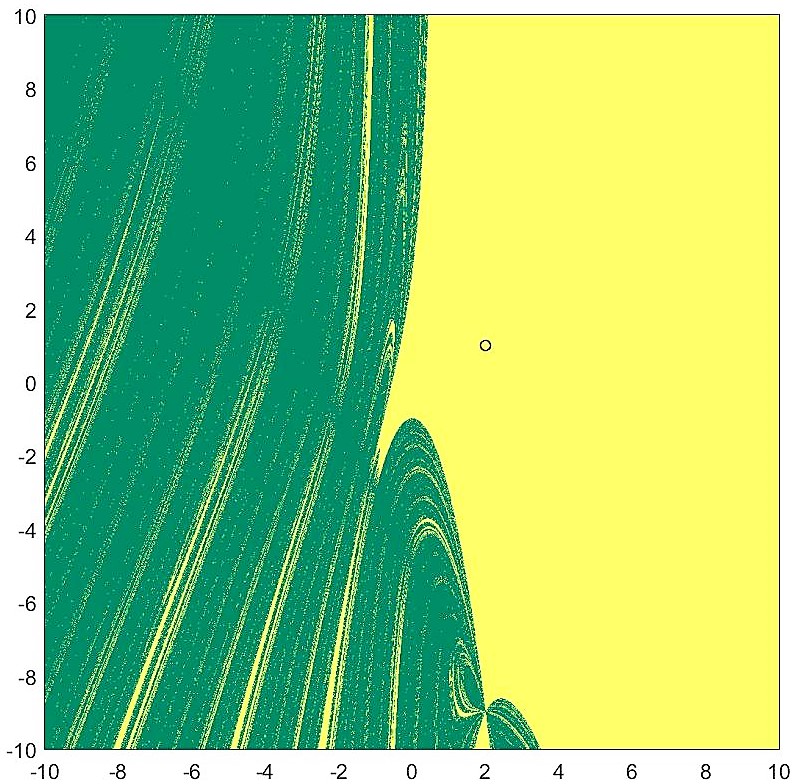}
\hfill
\includegraphics[width=0.45\textwidth]{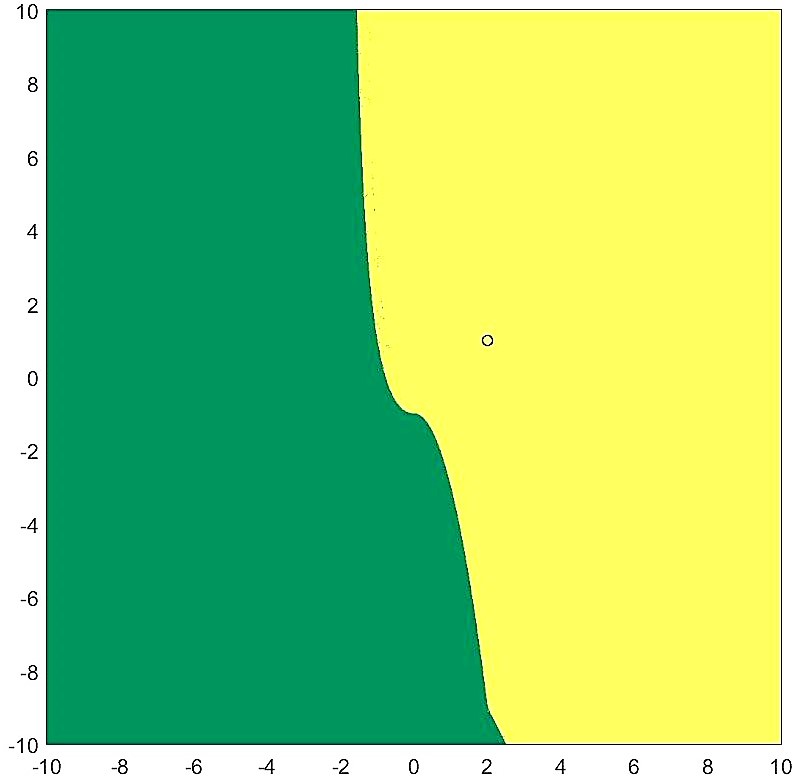}
\caption{The basins of attraction for Example~\ref{ex:3} by the Newton method. On the left without step size control (i.e.,~$t=1$) and on the right with step size control~($\tau = 0.01$). Note that the right part and yellow shaded domain indicates the domain, where the iteration procedure~\ref{al:full} converges to the unique root~$(2,1)$.}
\label{fig:frac56}
\end{figure}

\end{example}

\subsection{Performance data:}

\begin{figure}
\includegraphics[width=0.45\textwidth]{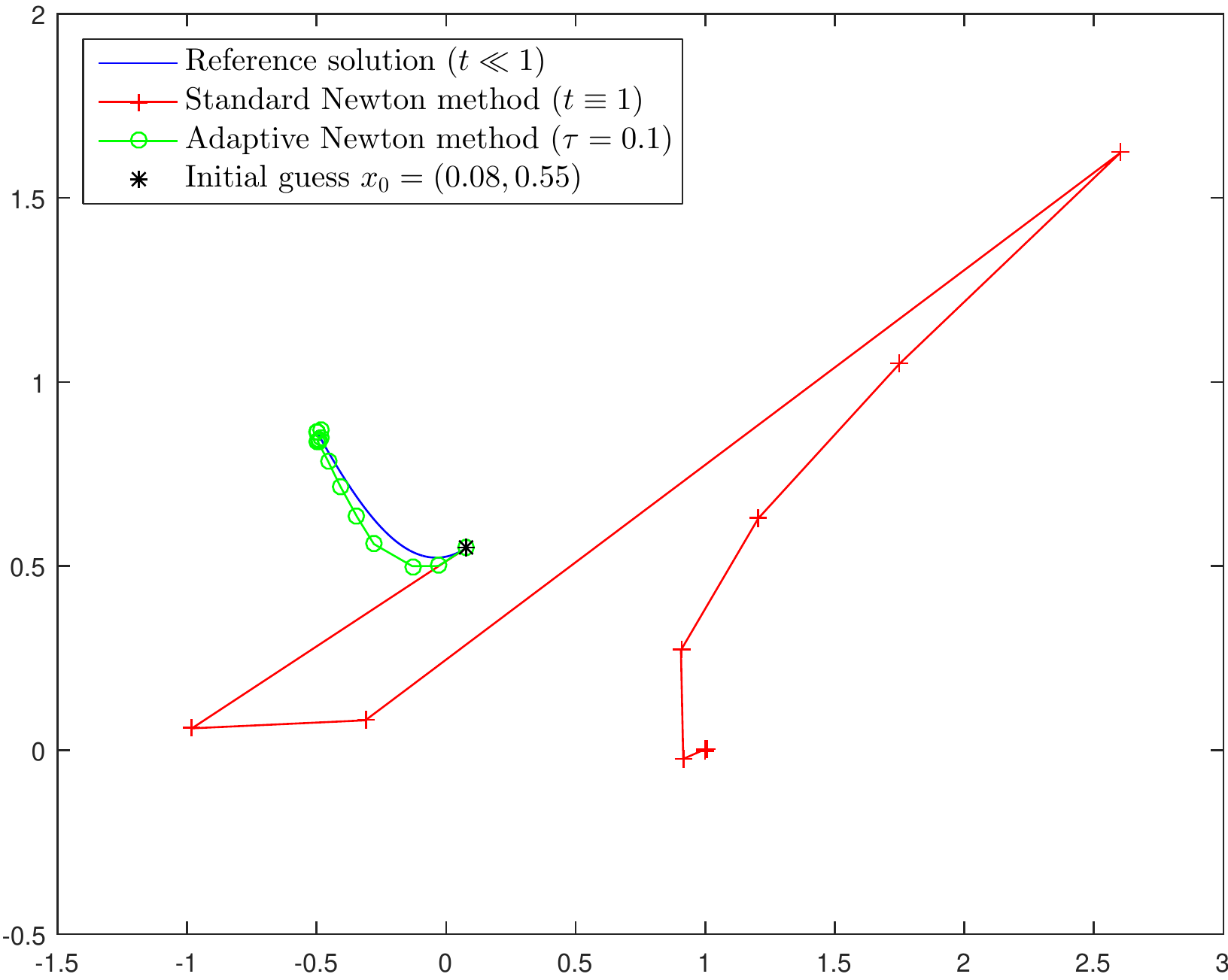}
\hfill
\includegraphics[width=0.45\textwidth]{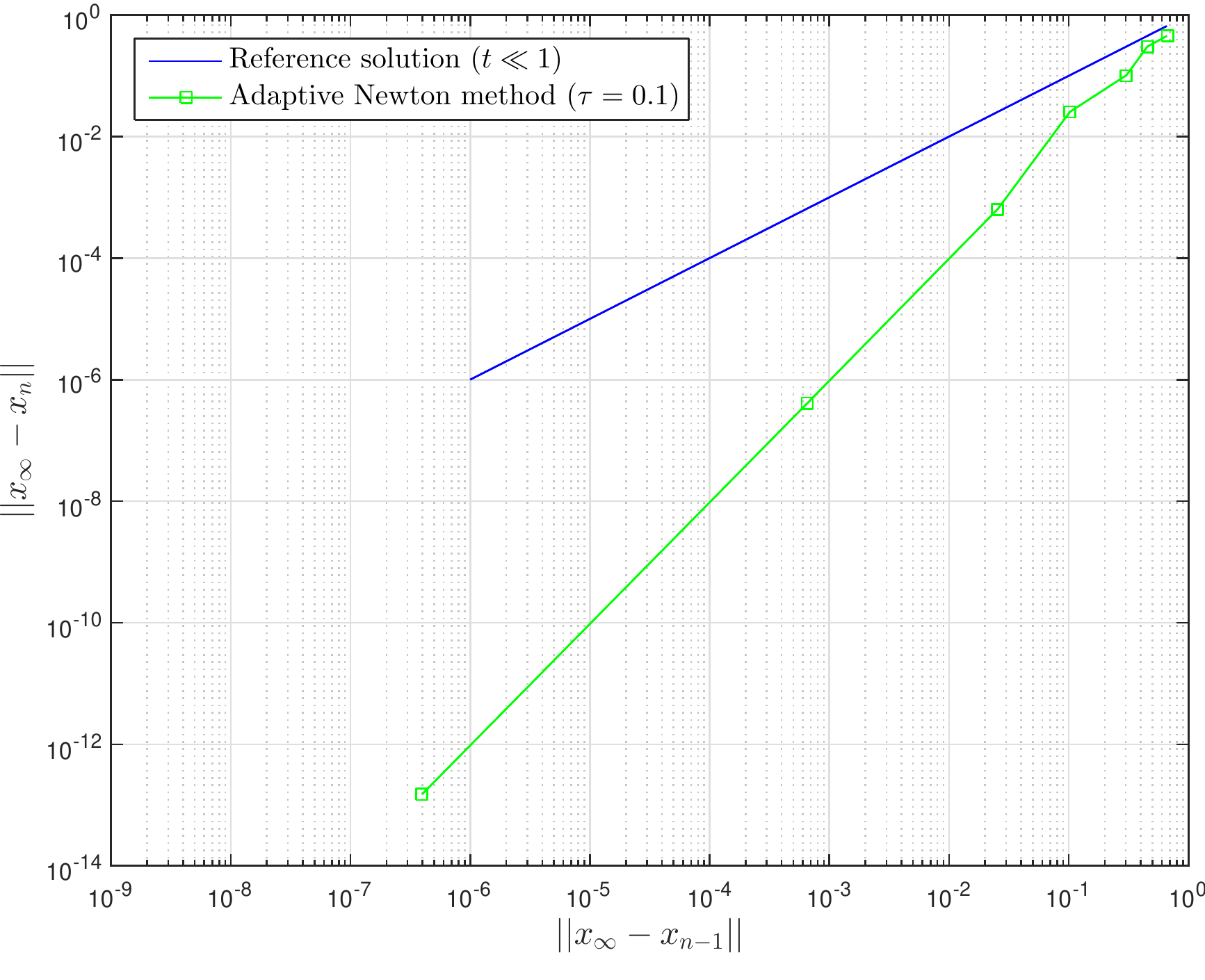}
\caption{Classical versus adaptive Newton method (left) and the convergence graphs corresponding to the reference solution and the adaptive iteration scheme (right).}
\label{fig:perf}
\end{figure}

\begin{figure}
\includegraphics[width=0.45\textwidth]{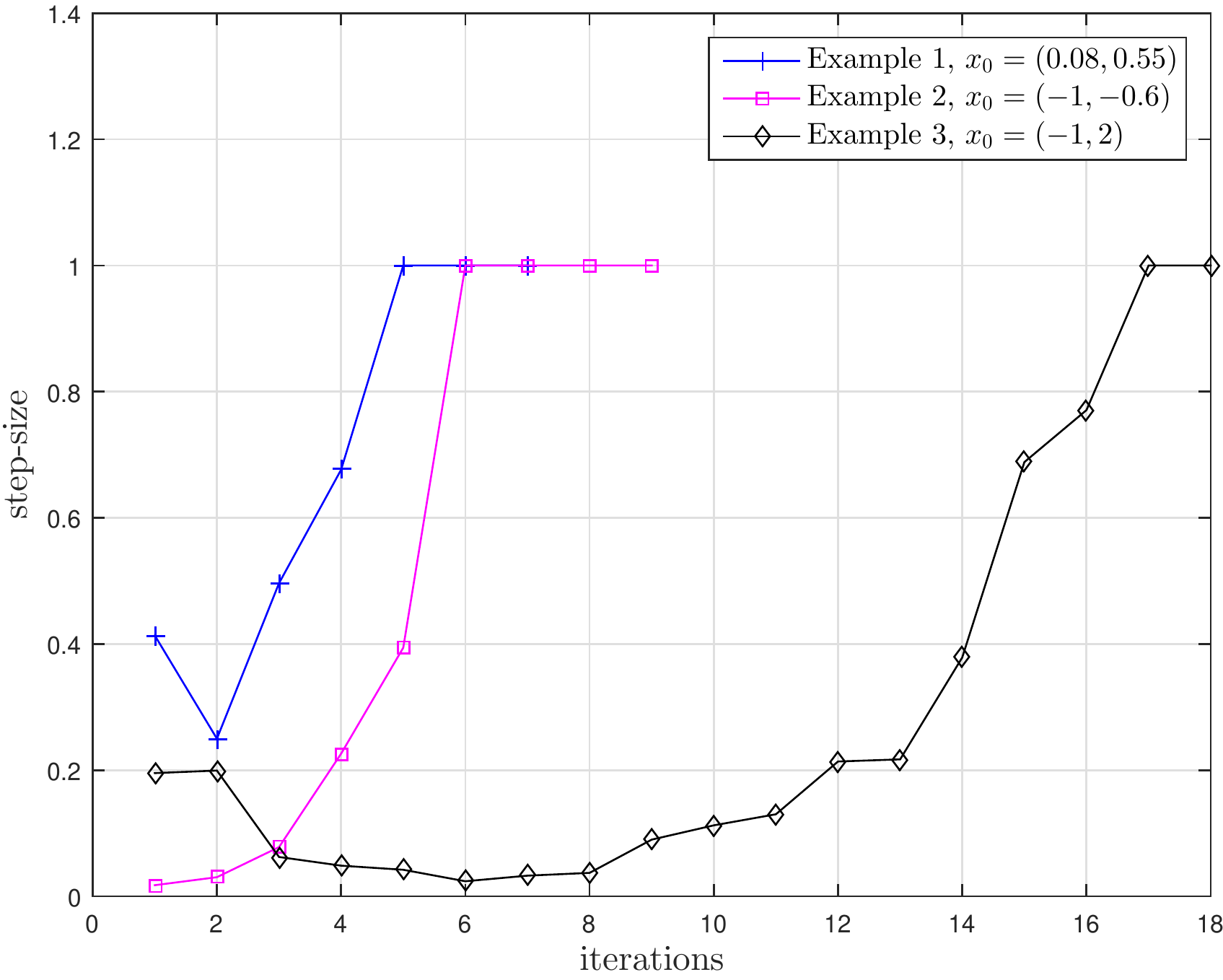}
\caption{Step-size versus number of effective computed updates in Algorithm~\ref{al:full} (Step 23). Here, $x_0$ denotes the initial value used in the depicted iteration (with $\tau=0.1$ and $\varepsilon = 10^{-9}$).}
\label{fig:step-size}
\end{figure}

\begin{table}
\begin{tabular}{p{5cm} l l c}
  & Example 1 on $[-3,3]^2$ & Example 2 on $ [-1.5,1.5]^2$ \\
\hline
 \% of convergent iterations with $t\equiv 1$. & 88.7\% & $55.44\% $ \\
\hline 
 \% of convergent iterations with adaptive step size $t$.  & $99.99\%$ & $70.5\%$ \\
\hline
\end{tabular}
\vspace{0.5cm}
\caption{Performance for Examples ~\ref{ex:1} and ~\ref{ex:2}.}
\label{tab:1}
\end{table}

In Figure~\ref{fig:perf}, we display the behavior of the classical and the adaptive Newton scheme (with~$\tau=0.1$). For Example~\ref{ex:2} we start the iteration in~$x_0=(0.08,0.55)$. Note that~$x_0$ is located in the---exact---attractor of the zero~$(-\nicefrac{1}{2},\nicefrac{\sqrt{3}}{2})$. We see that the classical solution with step size~$t \equiv 1$ shows large updates and thereby leaves the original attractor. On the other hand, the iterates generated by the adaptive scheme follow the exact solution (which was approximated by a numerical reference solution by choosing~$t\ll 1$) quite closely and is therefore able to approach the zero which is located in the corresponding domain of attraction.

On the right of Figure~\ref{fig:perf}, we show the convergence graphs corresponding to Example~\ref{ex:1} with the initial guess~$x_{0}=(0.08,0.55)$. Evidently, the adaptive iteration scheme shows quadratic convergence while the Newton scheme with fixed step size~$t\ll 1$ converges only linearly.

In Table~\ref{tab:1}, we depict the benefit of the proposed adaptive approach for Examples~\ref{ex:1} and~\ref{ex:2}. The numerical results in Table~\ref{tab:1} are based on the following considerations: For Examples~\ref{ex:1} and~\ref{ex:2}, we sample~$25 \times 10^4$ (equally-distributed) initial guesses on the domain~$[-3,3]^2$ and~$2.5\times 10^4$ on the domain~$[-1.5,1.5]^2$ respectively. Moreover, we call an initial value~$x_0$ convergent if it is in fact convergent and, additionally, approaches the correct zero, i.e. the zero that is  located in the same exact attractor as the initial value~$x_0$. The results in Table~\ref{tab:1} clearly show that the proposed adaptive strategy is able to enlarge the domain of convergence considerably.

Finally, let us address again Example~\ref{ex:3} in Table~\ref{tab:2}. These results are based on the following: Here, we call an initial value~$x_{0}$ convergent if it is in fact convergent, i.e., we skip the necessity that~$x_0$ belongs to the attractor of the unique root~$x_{\infty}=(2,1)$. This implies that the classical Newton method is now considered as convergent in subdomains of~$[-10,10]^2$ where the adaptive scheme is possibly not convergent (since for such an initial guess the trajectory of the continuous solution does not end at~$x_{\infty}$). Here, we sample~$10^6$ initial guesses on the domain~$[-10,10]^2$. In Table~\ref{tab:2}, we clearly see that the classical Newton method with step size~$t\equiv 1$ is convergent in~$51.2\%$ of the tested values while the adaptive scheme is only convergent in~$50.2\%$ of all cases. This fact nicely demonstrates that in certain situations a chaotic behavior of the iteration process is preferable in the sense that the iterates generated by the classical scheme are possibly able to cross critical interfaces with singular Jacobian. However, ---unnoticed---crossings between different basins of attraction and therefore a switching between different solutions of nonlinear problems can be considerably reduced by the proposed adaptive scheme.

\begin{table}
\begin{tabular}{p{5cm} l l}
& Example 3 on~$[-10,10]^2 $\\
\hline
\% of convergent iterations with $t\equiv 1$. & 51.2\%  \\
\hline 
\% of convergent iterations with adaptive step size~$t$.  & $50.2\%$ \\
\hline
\end{tabular}
\vspace{0.5cm}
\caption{Performance for Example~\ref{ex:3}.}
\label{tab:2}
\end{table}

\section{Conclusions}
\label{sec:concl}
In this work, we have considered an adaptive method for Newton iteration schemes for nonlinear equations,~$f(x)=0$ in~$\mathbb{R}^{n}$. The  adaptivity presented can be interpreted as a projection of a single iteration step onto the discretized global flow generated by the dynamics of the initial value problem~$\dot{x}=\A(x)f(x)$. Indeed, this system can be understood as a preconditioned version of the system~$\dot{x}=f(x)$ by~$\A(x)$. In particular, an appropriate choice of the matrix~$\A(x)$---if possible---can lead to the favorable property of all zeros being---at least on a local level---attractive. On the other hand---especially in case of~$\A(x)=\J_{f}(x)^{-1} $---singularities in~$\J_{f}$ may cause the associated discrete version to exhibit chaotic behavior. In order to tame these effects, we have used an adaptive step size control procedure whose purpose is to follow the flow of the \emph{continuous} system to a certain extent. We have tested our method on a few low dimensional problems. Moreover, our experiments demonstrate empirically that the proposed scheme is indeed capable to tame the \emph{chaotic} behavior of the iteration compared with the classical Newton scheme, i.e., without applying any step size control. In particular, our test examples illustrate that high convergence rates can be retained, and the domains of convergence can---typically---be considerably enlarged.

\bibliographystyle{amsplain}
\bibliography{references}
\end{document}